\documentclass[12pt]{article}
\usepackage{graphicx}
\usepackage{amsmath,amssymb,amsfonts,latexsym,amsthm}

\newcommand{\virg}{\ ,}
\newcommand{\pvirg}{\ ;}

\newcommand{\ponto}{\ .}
\newcommand{\timesd}{{\times}}
\newcommand{\N}{\mathbb{N}}

\theoremstyle{definition}
\newtheorem{exmp}{Example}[section]

\begin{document}

\title{Fermat's Last Theorem admits an infinity of proving ways and two
corollaries}
\author{Jos\'{e} Cayolla \\
(Jubilated Professor, Instituto Superior de Engenharia de Lisboa)\\
Rua de S. Paulo, 107, 1%
${{}^o}$%
F, 2775-752 Carcavelos, Portugal\\
Tel.: +351-214539661\\
email: josecayolla@sapo.pt}
\date{ }
\maketitle

\begin{abstract}
Fermat's statement is equivalent to say that if $x$, $y$, $z$, $n$ are
integers and $n>2$, then $z^{n}\gtrless x^{n}+y^{n}$. This is proved with
the aid of numbers $\lambda $'s, of the form $\lambda =z/\rho $, with $%
1<\rho <z$, named \emph{reversors} in the text, because their property of
multiplying $z^{n-1}$ in $z^{n-1}<x^{n-1}+y^{n-1}$, not only reverses the
signal but also gives $z^{n}>x^{n}+y^{n}$ as a solution of the reversed
inequality. As the $\lambda ^{\prime }s$ satisfy a compatible opposed sense
system of inequalities, the $\lambda $-set is equivalent to the points of an
$\mathbb{R}^{+}$ interval. THerefore the theorem admits a noncountable
infinity of proving ways, each one given by a particular value of $\lambda $.

In Corollary 1 a general relation between $y$, $x$, $z$ and $n$ is derived.
Corollary 2 shows that the Diophantine equation in Fermat's statement admits
no solutions other than algebraic irrationals and the inherent complexes.

Integer triplets can be classified in seven sets, within each one their
relation with the respective $n$ is the same as shown in Table~\ref{tab:1}.

Numerical verification with examples taken from all the mentioned seven sets
gives a total agreement with the theory.
\end{abstract}

\section{Introduction}

\label{intro}

The history of Fermat's Last Theorem (FLT) is so well known, that we
restrict ourselves to remind only some important historical facts.

In 1637 Fermat wrote his famous marginal note in Diophanto's ``Arithmetica",
stating the theorem together with the comment that he ``had found a
remarkable proof, but the margin's space was not enough to write it there"
\cite[p.\,173]{b3}.

Strangely enough, Fermat let elapse 28 years without publishing any proof
until he died in 1665 \cite[p.\,353]{b3}.

His son published the note in 1670 \cite[p.\,172]{b3}, so starting the FLT
saga.

The theorem resisted the efforts of the best mathematicians for a lot more
than three hundred years, so much so, that in the seventies of last century
it was wide spread the opinion that Fermat realized that what he had in mind
was a wrong or an incomplete proof \cite[p.\,202]{b2}.

Finally, in the late nineties, the british mathematician Andrew Wiles
achieved a full FLT proof. A very brilliant and sophisticated one. Indeed,
the way he chose was through a previous proof of Shimura--Taniyama--Weil
conjecture, not yet achieved at that time, considered more general than FLT,
even more difficult and implying the FLT proof. This approach gave Wiles the
opportunity of producing a lot of remarkable contributions mainly in the
fields of elliptic curves modularity, topology and number theory \cite{b1},
\cite{b4}. Obviously, this asked for an extensive use of outstanding modern
mathematical concepts, all this work deserving congratulations, because the
resulting progress of nowadays mathematics. Fermat did a similar thing
regarding the mathematics of his time, but failed in his ``Last Theorem".

Still an interesting historical note is that \emph{FLT in the special case $%
n=3$, was already known in ancient Greece} under the name of ``\emph{the
cube's duplication problem}". \emph{Mathematically}, the rising of this
question is quite natural. Indeed, Pythagoras, with his fundamental theorem,
had \emph{duplicated the square}. So, why not to start thinking of
duplicating a cube? In other words, given an integer and its cube, to find
two other integers such that the sum of their cubes equals the first cube.
Notorious intellectuals, such as Arquitas, Democritus and many others,
became interested, but none solved the problem. \emph{The situation remained
so for more than two thousand years}. Then, a century after Fermat's
marginal note, \emph{Euler has proved that the problem was an impossible one}%
.

Reminding all we said now, we kept the opinion that it seems quite
defensible to think that achieving a full and very simple FLT proof still
remains an open and quite challenging problem.

After some time and unsuccessful attempts, we have been lucky enough to find
this way.

\section{Statement}

\label{sec:1}

The equation $z^n=x^n+y^n$ admits no integer solutions in $x,y$ and $z$ if $n
$, also an integer, is greater than 2. This is equivalent to the following.
\begin{equation}  \label{Eq1.1}
\hbox{If $\{y,x,z,n\}\in\N^{+}$ and $n>2$, then $z^{n}\gtrless
x^{n}+y^{n}$\ponto}
\end{equation}

\section{Preliminary considerations}

\label{sec:2}

From now on and just for short, we refer to the theorem writing FLT.

By definition, we say that the inequality between the $z$ power and the sum
of the other two, with the same exponent, is \emph{unreversed} if the signal
is $<$ and \emph{reversed} if $>$.

When it happens $z^{n-1}<x^{n-1}+y^{n-1}$ and $z^n>x^n+y^n$, we call $n$ the
\emph{reversion exponent}. Supposing that there are solutions of the
Diophatine equation form in the statement, its form \eqref{Eq1.1}, implies $%
z>x$ and $z>y$.

In conformity with these implications are the two hypothesis $z\geqslant x+y$%
. However, $z>x+y$ remains reversed with $n\geqslant 1$. The same happens to
$z=x+y$ with $n>1$.

The complementar of $z\geqslant x+y$ is $z<x+y$. Now the representative
segments of $x,y$ and $z$ may define a triangle. Triangles can be classified
through the parameter \emph{greatest inner angle, $\alpha$}. If $\alpha>\pi/2
$, $z$ is too large to define a Pythagorean triplet with $x$ and $y$. Hence,
$z^2>x^2+y^2$, so $n=2$.

With $\alpha=\pi/2$ we have a right triangle. Its inequality reverses
through an intermediary equality, $z^2=x^2+y^2$. But, from $zx^2>x^3$ and $%
zy^2>y^3$ we have $z^3>x^3+y^3$, so $n=3$. We remark that $x=y$ \emph{gives $%
z$ irrational, not integer}, with exponents one and any other odd number.

It follows $\alpha<\pi/2$.

This subdivides in three cases. One with all inner angles $\alpha=\pi/3$, so
$z=x=y$, contradicting $z>x$ and $z>y$, then no solutions here.
Geometrically, the inequality $z<z+z$ never reverses, whatever the exponent,
thus, also no $n$. Another similar situation happens with two $\alpha$ and
the opposite sides $z$ equal, contradicting $z>x$ and $z>y$. Then again no
solutions and no $n$.

From what have been said and also from the symbology about $x$ and $y$ being
arbitrary, we can assume with no loss of generality,
\begin{equation}  \label{eq2.1}
z>x>y\ .
\end{equation}
Therefore, the only remaining triangle's type is defined by $\alpha<\pi/2$
and all sides different. Now $z$ is not large enough to define a Pythagorean
triplet with $x$ and $y$. Then, the inequations defining this triangle are
\begin{equation}  \label{eq2.2}
\begin{aligned} z&< x+y\ ,\\ z^{2} &< x^{2}+y^{2}\ , \end{aligned}
\end{equation}
the second showing that \emph{to get out of the situation unreversed}, an
exponent, \emph{depending on $y,x$ and $z$ and greater than two} is needed,
\begin{equation}  \label{eq2.3}
n(y,x,z)>2\ .
\end{equation}
What have been said up to now is synthesized in Table \ref{tab:1} below, the
symbol $\triangle$ meaning ``triangle".

\begin{table}[h]
\caption{Triplets' classification.}
\label{tab:1}
\begin{tabular}{|llllllllllll|}
\hline
\multicolumn{6}{|l|}{Set 1} & \multicolumn{6}{|l|}{Set 2} \\
\multicolumn{6}{|l|}{$\nexists\,\triangle\,\leftrightarrow\,z\ge x+y$} &
\multicolumn{6}{|l|}{$\exists\,\triangle\,\leftrightarrow\,z<x+y$} \\ \hline
\multicolumn{3}{|l|}{Subset 1.1} & \multicolumn{3}{|l|}{Subset 1.2} &
\multicolumn{3}{|l|}{Subset 2.1} & \multicolumn{3}{|l|}{Subset 2.2} \\
\multicolumn{3}{|l|}{$z>x+y$} & \multicolumn{3}{|l|}{$z=x+y$} &
\multicolumn{3}{|l|}{$\alpha>\pi/2$} & \multicolumn{3}{|l|}{$\alpha=\pi/2$}
\\
\multicolumn{3}{|l|}{$n=1$} & \multicolumn{3}{|l|}{$z^{2}>x^{2}+y^{2}$} &
\multicolumn{3}{|l|}{$z^{2}>x^{2}+y^{2}$} & \multicolumn{3}{|l|}{$%
z^{2}=x^{2}+y^{2}$} \\
\multicolumn{3}{|l|}{} & \multicolumn{3}{|l|}{$n=2$} & \multicolumn{3}{|l|}{$%
n=2$} & \multicolumn{3}{|l|}{$z^{3}>x^{3}+y^{3}$} \\
\multicolumn{3}{|l|}{} & \multicolumn{3}{|l|}{} & \multicolumn{3}{|l|}{} &
\multicolumn{3}{|l|}{$n=3$} \\ \hline
\multicolumn{4}{|l|}{Subset 2.3} & \multicolumn{4}{|l|}{Subset 2.3.1} &
\multicolumn{4}{|l|}{Subset 2.3.2} \\
\multicolumn{4}{|l|}{$\alpha<\pi/2$} & \multicolumn{4}{|l|}{$\alpha<\pi/2$;
\,$z=x>y$} & \multicolumn{4}{|l|}{$\alpha<\pi/2$; \,$z=x=y$} \\ \hline
\multicolumn{4}{|l|}{Subset 2.3.1} & \multicolumn{4}{|c|}{$\nexists\,n$,} &
\multicolumn{4}{|c|}{$\nexists\,n$,} \\
\multicolumn{4}{|l|}{$z^{2}<x^{2}+y^{2}$, \,$n\ge 3$ and} &
\multicolumn{4}{|l|}{$\triangle$ type conserved} & \multicolumn{4}{|l|}{$%
\triangle$ type and angles} \\
\multicolumn{4}{|l|}{depending on $\{y,x,z\}$} & \multicolumn{4}{|l|}{but
not angles} & \multicolumn{4}{|l|}{conserved} \\ \hline
\end{tabular}%
\end{table}

It is now clear that proving FLT in general terms is reduced to do that in
the special terms \eqref{eq2.2}.

\section{Proof}

\label{sec:3}

\subsection{A necessary and sufficient condition}

\label{sec:3.1}

One of such conditions is the existence of two consecutive integers, $n-1$
and $n$, \emph{depending on $y,x$ and $z$}, and such that
\begin{equation}
\begin{aligned} z^{n-1}&< x^{n-1}+y^{n-1}\ ,\\ z^{n} &> x^{n}+y^{n}\ .
\end{aligned}  \label{eq3.1}
\end{equation}%
An alternative and equivalent way of stating \eqref{eq3.1} is the following.
If we imagine the exponent as being the continuous variable in $\mathbb{R}%
^{+}$, the implicit function
\begin{equation}
z^{s}-x^{s}-y^{s}=0  \label{eq3.2}
\end{equation}%
defines $s$ in terms of $y,x$ and $z$ and with the requirement
\begin{equation}
s\in ]n-1,n[  \label{eq3.3}
\end{equation}

\subsection{Accessory entities}

\label{sec:3.2}

For $\forall\, i\in N^{+}$, we define
\begin{equation}  \label{eq3.4}
p_i=x^{i}+y^{i}\ .
\end{equation}
From \eqref{eq2.1} $x>y$, then
\begin{equation}  \label{eq3.5}
\begin{aligned} p_{i}\,x &> p_{i+1}\ ,\\ p_{i}\,y &< p_{i+1}\ , \end{aligned}
\end{equation}
so,
\begin{equation}  \label{eq3.6}
\frac{p_{i+1}}{p_i}=\frac{x^{i+1}+y^{i+1}}{x^{i}+y^{i}}=k_{i}\ .
\end{equation}
This, together with \eqref{eq3.5} show that the $k_{i}$ are rational
numbers, increasing with $i$, but always satisfying
\begin{equation}  \label{eq3.7}
y<k_{i}<x\ .
\end{equation}
We can now compare the successions of the $z^{i}$ powers with that of the $%
p_{i}$,
\begin{equation}  \label{eq3.8}
\begin{aligned} z^{i+1}&=&z^{i}z\\ p_{i+1}&=&p_{i}k_{i} \end{aligned}
\end{equation}
showing that the $z^{i}$ grow ``faster" than the $p_{i}$. Thus, there is a
certain $i$ such that $z^{i}>p_{i}$, and so for any further $i$. By other
words, \emph{when the exponents increase}, $z<x+y$ in \eqref{eq2.2} reverses
and remains so for any further exponent. This allows to define the number $%
n-1$ in \eqref{eq3.1} as follows:
\begin{equation}  \label{eq3.9}
n-1=\max\{i\in N^{+}: z^{i}<x^{i}+y^{i}\}
\end{equation}

\subsection{Geometrical evolution}

\label{sec:3.3}

When the exponents increase the initial triangle \eqref{eq2.2} changes size
and shape up to exponents $n-1$, with which there is still a triangle whose
inequality is
\begin{equation}  \label{eq3.10}
z^{n-1}<x^{n-1}+y^{n-1}\ .
\end{equation}
With exponents $n$, \emph{be FLT true or false}, there is no triangle,
because $z^n=x^n+y^n$, or $z^n>x^n+y^n$, corresponding the first to FLT
false, and the second to FLT true, \emph{define no triangle}. And, as
already shown by \eqref{eq3.8} and \eqref{eq3.9}, the inequality %
\eqref{eq3.10} is reversed \emph{with any exponent greater than $n$}.

This gives two conclusions:

a) without contradiction, FLT can only be supposed false with exponents $n$.

b) with exponents $n-1$ the triangle
\begin{equation}  \label{eq3.11}
T_{n-1}=\{y^{n-1},x^{n-1}, z^{n-1}\}
\end{equation}
can be called ``the last triangle". Its nature may be found squaring the
sides, giving exponents $2n-2$. To the inequality $z<x+y$ being surely
reversed it must be $2n-2>n$, this asking for $n>2$, which is precisely the
case as shown by \eqref{eq2.2} and \eqref{eq2.3}. Then,
\begin{equation}  \label{eq3.12}
z^{2n-2}>x^{2n-2}+y^{2n-2}\ ,
\end{equation}
proving that in the last triangle, $T_{n-1}$ of \eqref{eq3.11}, $\alpha>\pi/2
$; \emph{independent of $y$, $x$ and $z$}.

We remark that the reversion in the squares, \eqref{eq3.12}, is not
``strong" enough to prove FLT, because this theorem \emph{asks for a
reversion of \eqref{eq3.10} in the next integer exponent $n$}.

\subsection{Reversors}

\label{sec:3.4}

About the necessary and sufficient condition \eqref{eq3.1}, we have already
proved the ever existence of the number $n-1$, defined by \eqref{eq3.9}.
Hence, it remains to prove the second of \eqref{eq3.1}, which can be written
in the form
\begin{equation}  \label{eq3.13}
z^n-x^n-y^n>0\ .
\end{equation}
This positive difference, \emph{if FLT is true}, defines an interval in $%
\mathbb{R}^{+}$, of numbers, $\zeta_n$. Thus, we can write
\begin{equation}  \label{eq3.14}
\text{FLT true} \Longleftrightarrow \exists \{\zeta_n\}\,:\,z^n\geqslant
\zeta_n\geqslant x^n+y^n\ .
\end{equation}
This clearly shows that ``FLT is true" is a proposition \emph{admiting a
continuous infinity of proving ways}, each one through a particular value of
$\zeta_n$. Hence, the next step is to find an algorithm linking $y$, $x$ and
$z$ with the $\zeta_n$. We start to do that searching for a continuous set
of numbers, $\lambda^{\lambda}$, which we call \emph{reversors}, because of
their property of, \emph{despite smaller than $z$}, multiplying by a $\lambda
$ the $z^{n-1}$ in the last triangle's inequality, \eqref{eq3.10}, \emph{%
contradicts the situation unreversed, and, at the same time, gives the
solution $z^n>x^n+y^n$, to the new situation}.

Formally, as the contradiction of $<$ is $\geqslant$, this property is
written
\begin{equation}  \label{eq3.15}
\lambda\,z^{n-1}\geqslant x^{n-1}+y^{n-1}\Longrightarrow z^n>x^n+y^n\ .
\end{equation}
It is clear that from \eqref{eq3.15}, \emph{to be proved soon}, the
existence of \emph{at least one $\lambda$ is a sufficient condition} of FLT
being true. To have the property \eqref{eq3.15} a $\lambda$ must be of the
form
\begin{equation}  \label{eq3.16}
1<\lambda=\frac{z}{\rho}<z\ ,
\end{equation}
implying
\begin{equation}  \label{eq3.17}
1<\rho<z\ .
\end{equation}
Entering with this in \eqref{eq3.15}, we have
\begin{equation}  \label{eq3.18}
\frac{z}{\rho} \, z^{n-1}\geqslant x^{n-1}+y^{n-1}\Longrightarrow
z^n>x^n+y^n\ ,
\end{equation}
or
\begin{equation}  \label{eq3.19}
\rho(x^{n-1}+y^{n-1})\geqslant x^n+y^n\ ,
\end{equation}
giving
\begin{equation}  \label{eq3.20}
\rho\geqslant\frac{x^n+y^n}{x^{n-1}+y^{n-1}}=\frac{p_n}{p_{n-1}}=k_{n-1}\ ,
\end{equation}
having retaken the use of the $p_i$ and $k_i$, defined in Sec. \ref{sec:3.2}%
. The $kk$ property \eqref{eq3.7}, $y<k_i<x<z$, confirms that $k_{n-1}$ is a
$\rho$. More precisely, $\rho\geqslant k_{n-1}$ shows that
\begin{equation}  \label{eq3.21}
k_{n-1}=\min\{\rho\}=\underline{\rho}\ .
\end{equation}
by definition of $\underline{\rho}$.

Hence,
\begin{equation}
\exists \,\lambda =\frac{z}{k_{n-1}}\ ,  \label{eq3.22}
\end{equation}%
\emph{sufficient to prove that FLT is true}. Indeed, entering with $\rho
=k_{n-1}$ in \eqref{eq3.18} and performing the calculations gives
\begin{equation}
z^{n}>x^{n}+y^{n}\ ,  \label{eq3.23}
\end{equation}%
\emph{like we intended to prove}. As $k_{n-1}=\underline{\rho }$,
\begin{equation}
\frac{z}{k_{n-1}}=\max \{\lambda \}\ ,  \label{eq3.24}
\end{equation}%
and this rises the question of what could be $\max \{\rho \}=\overline{\rho }
$.

Dividing \eqref{eq3.18} by $z^{n-1}$, we obtain
\begin{equation}  \label{eq3.25}
z>\lambda\geqslant\frac{x^{n-1}+y^{n-1}}{z^{n-1}}=\varphi>1\ ,
\end{equation}
having introduce the symbol $\varphi$, which is greater than one by force of %
\eqref{eq3.10}. The signal $\geqslant$ in \eqref{eq3.25} gives
\begin{equation}  \label{eq3.26}
\min\{\lambda\}=\varphi=\frac{x^{n-1}+y^{n-1}}{z^{n-1}}\ ,
\end{equation}
then,
\begin{equation}  \label{eq3.27}
\varphi=\frac{z}{\overline{\rho}}
\end{equation}
and
\begin{equation}  \label{eq3.28}
\overline{\rho}=\frac{z}{\varphi}=\frac{z^n}{x^{n-1}+y^{n-1}}\ .
\end{equation}

To verify \eqref{eq3.26}, we enter with \eqref{eq3.28} in \eqref{eq3.25},
giving
\begin{equation}  \label{eq3.29}
z^n=\overline{\rho}(x^{n-1}+y^{n-1})=z^n\ ,
\end{equation}
confirming, not only that $\varphi=\min\{\lambda\}$, but also that $\varphi$
reverses \eqref{eq3.10} \emph{in the strongest possible way}, because the
reversion's result is $z^n$ itself, the $\max\{\zeta_n\}$. All these
considerations are verifiable by the relations
\begin{equation}  \label{eq3.30}
\overline{\rho}>\underline{\rho}
\end{equation}
or
\begin{equation}  \label{eq3.31}
\frac{z}{k_{n-1}}>\varphi\ .
\end{equation}
Indeed, attending the $k_{n-1}$ form, in \eqref{eq3.20}, supposing $%
z/k_{n-1}<\varphi$ gives
\begin{equation}  \label{eq3.32}
z^n<x^n+y^n\ ,
\end{equation}
contradicting \eqref{eq3.9}.

The alternative, $z/k_{n-1}=\varphi$ gives
\begin{equation}  \label{eq3.33}
z^n=x^n+y^n\ ,
\end{equation}
i.e. \emph{FLT false, with exponents $n$}, contradicting FLT true, as
already proved by \eqref{eq3.22}. Hence, it only remains the conclusion that
$z/k_{n-1}>\varphi$, as we wanted to verify. Therefore, the main conclusion
is that FLT \emph{is also true} for any triplet of the type \eqref{eq2.2}.

As in Sec. \ref{sec:2} above, FLT has been proved in all other triplets
shown in Table \ref{tab:1}, \emph{this FLT full proof is now completed}.

Another main conclusion \emph{is that FLT admits a not countable infinity of
proving ways}, each one through a particular value of the number $\lambda $
belonging to the set
\begin{equation}
\{\lambda \}\Longleftrightarrow \left[ \varphi ,\frac{z}{k_{n-1}}\right] \ .
\end{equation}%
\emph{Another way showing the existence of a continuous infinity of proving
ways is the fact that $\lambda ^{\prime }s$ are the solutions of the
compatible opposite sense inequations system \eqref{eq3.25}}.

\subsection{Overreversors}

\label{sec:3.5}

The algorithm linking $y$, $x$ and $z$ with the $\rho ^{\prime }s$, $\lambda
^{\prime }s$, and $\zeta ^{\prime }s$ is the following.

Choosing a $\rho$ from the set
\begin{equation}  \label{eq3.35}
\{\rho\}\Longleftrightarrow \left[k_{n-1},\frac{z}{\varphi}\right]\ ,
\end{equation}
we have the corresponding $\lambda=z/\rho$ and $\zeta_n=\rho(x^{n-1}+y^{n-1})
$.

Double opening the $\rho $ interval, \eqref{eq3.35}, the same happens with
those of the $\lambda ^{\prime }s$ and $\zeta ^{\prime }s$. Therefore, to
\begin{equation}
\forall \,\lambda \in \left] \varphi ,\frac{z}{k_{n-1}}\right[
\label{eq3.36}
\end{equation}%
corresponds a $\zeta _{n}$ such that
\begin{equation}
z^{n}>\zeta _{n}>x^{n}+y^{n}\ ,  \label{eq3.37}
\end{equation}%
\emph{logically stronger than} just $z^{n}>x^{n}+y^{n}$.

Because of this property, we call the $\lambda ^{\prime }s$ from %
\eqref{eq3.36} \emph{overreversors} and the respective $\zeta _{n}$ \emph{an
overreversion}, because more than needed to prove FLT. For example, if $z$
is large enough as $z^{n}>x^{n}+xy^{n-1}$, $z/x$ is an overreversor.

A numerical example of this situation is the following
\begin{align*}
&\{y,x,z\}=\{2,3,4\}; \\
&4<3+2; \;4^2>3^2+2^2; \\
&n=2;\; n-1=1;\; T_{n-1}=T_1=\{2,3,4\}; \\
&\varphi=\frac{3+2}{4}=\frac{5}{4}; \\
&k_{n-1}=k_1=\frac{3^2+2^2}{3+2}=\frac{13}{5}; \\
&\frac{z}{k_1}=\frac{20}{13}; \\
&\{\lambda\}=\left[\frac{5}{4},\frac{20}{13}\right]; \\
&z^2>k_1(x+y);\; 4^2>\frac{13}{5}(3+2); \; 16>13; \\
&\frac{z}{x}=\frac{4}{3}\in\left]\frac{5}{4}\, \frac{20}{13}\right[; \\
&z^2>x^2+xy; \\
&\zeta_2=x(x+y);\; \zeta_2=15; \\
&4^2>\zeta_2>3^2+2^2; \; 16>15>13.
\end{align*}

\emph{Numerical examples, covering all the possibilities shown in Table \ref%
{tab:1}, are given in Sec. \ref{sec:6} below, in order to verify the
theoretical results achieved in this work}.

\section{Corollary 1, a general relation between $y,x,z$ and $n$.}

\label{sec:4-Corollary}

We have proved FLT through \eqref{eq3.1}. Taking Neperian logarithms and
solving in $n$ and $n-1$ gives
\begin{equation}  \label{eq5.1}
n>\frac{\log(x^n+y^n)}{\log z}=b
\end{equation}

\begin{equation}  \label{eq5.2}
n-1<\frac{\log(x^{n-1}+y^{n-1})}{\log z}=a\ .
\end{equation}
The first shows that through simple calculations we get $n$, and then $b$,
and again $n$, \emph{as the first integer above $b$}, so
\begin{equation}  \label{eq5.3}
0<n-\frac{\log(x^n+y^n)}{\log z}<1
\end{equation}
is the wanted relation.

To remark that $b$ is never an integer. However, it is possible $a=2$, if
and only if $y$, $x$ and $z$ define a Pythagorean triplet, as easily seen
through \eqref{eq5.2}.

From \eqref{eq5.1} $b<n$, then
\begin{equation}  \label{eq5.4}
z^b=x^n+y^n>x^b+y^b\ ,
\end{equation}
showing that \emph{the inequality \eqref{eq3.10} is already reversed with
exponents $b$}. In parallel terms, apart the Pythagorean triplets, $a>n-1$
and, from \eqref{eq5.2},
\begin{equation}  \label{eq5.5}
z^a=x^{n-1}+y^{n-1}<x^a+y^a\ ,
\end{equation}
proving that the inequality is not yet reversed with exponents $a$.

\subsection{Bounding $b-a$}

\label{sec:4.1}

Through \eqref{eq5.1} and \eqref{eq5.2} we clearly see that
\begin{equation}  \label{eq5.6}
0<b-a<1\ .
\end{equation}
Performing the subtraction \eqref{eq5.1}-\eqref{eq5.2} and reminding the $k_i
$ in \eqref{eq3.7}, one easily obtains
\begin{equation}  \label{eq5.7}
b-a=\frac{\log(k_{n-1})}{\log z}\ ,
\end{equation}
confirming \eqref{eq5.6}, because
\begin{equation}  \label{eq5.8}
y<k_{n-1}<x<z\ .
\end{equation}
By force of $0<b-a<1$, it is also
\begin{equation}  \label{eq5.9}
0<1-(b-a)<1\ .
\end{equation}

Attending \eqref{eq5.7}, we have
\begin{equation}  \label{eq5.10}
1-(b-a)=\frac{\log(z/k_{n-1})}{\log z}\ .
\end{equation}
This and \eqref{eq5.7} ask to compare $z/k_{n-1}$ with $k_{n-1}$. The best
way is through the $k_i$ property \eqref{eq3.7} and the inequations %
\eqref{eq2.2} of the initial triangle. This gives
\begin{equation}  \label{eq5.11}
k^2_{n-1}>y^2>(z+x)(z-x)\ .
\end{equation}

As $y,x$ and $z$ are integers, and, from \eqref{eq2.1} $z>x$, $z-x\geqslant 1
$, so $k^2_{n-1}>z$, or
\begin{equation}  \label{eq5.12}
\frac{z}{k_{n-1}}<k_{n-1}\ .
\end{equation}

Entering now with $z/k_{n-1}<k_{n-1}$ in \eqref{eq5.10} and comparing %
\eqref{eq5.12} with \eqref{eq5.7}, gives $1-(b-a)<b-a$, or
\begin{equation}  \label{eq5.13}
\frac{1}{2}<(b-a)<1\ ,
\end{equation}
the wanted better $(b-a)$ bounding.

By the same stroke this allows to improve \eqref{eq5.3}, because $n-b$ is a
part of $1-(b-a)$. So we can write \eqref{eq5.3} in the better form
\begin{equation}  \label{eq5.14}
0<n-\frac{\log(x^n+y^n)}{\log z}<\frac{1}{2}\ .
\end{equation}

\subsection{Bounding $s$}

In Sec. \ref{sec:3.2} we have introduced the equalizing exponent $s$,
implicitly defined by $z^s=x^s+y^s$. Taking Neperian logarithms and solving
in $s$, we have
\begin{equation}  \label{eq5.28}
s=\frac{\log (x^s+y^s)}{\log z}\ .
\end{equation}
By force of \eqref{eq5.5} and \eqref{eq5.6}, the $s$ bounds are
\begin{equation}  \label{eq5.29}
n-1\leq a<s<b<n\ .
\end{equation}
The main conclusions about this are synthesized in Fig. 1 below.
\begin{figure}[htb!]
\includegraphics[scale=0.84]{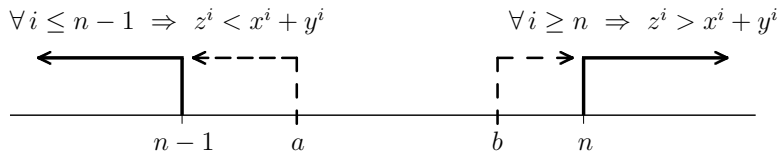} 
\caption{Inequality's distribution in $\mathbb{R}^{+}$. (No scales)}
\label{fig:1}
\end{figure}

\section{Corollary 2}

\label{sec:5-Corollary}

\emph{The Diophantine equation in Fermat's statement admits no solutions
other than algebraic irrationals and the inherent complexes}.

\subsection{Positive rationals}

\label{Sc7.1}

Supposing this corollary false and letting
\begin{equation}  \label{Eq7.1}
i=1,2,3;\ \ \ \{q_{i},q_{i}^{\prime }\}\in\mathbb{N}^{+}\ \ \hbox{and}\ \
n>2\ ,
\end{equation}
we could write
\begin{equation}  \label{Eq7.2}
(q_{3}/q_{3}^{\prime n}=(q_{2}/q_{2}^{\prime n}+(q_{1}/q_{1}^{\prime n}\ .
\end{equation}
Multiplying by $q_{1}^{\prime }q_{2}^{\prime }q_{3}^{\prime }$, we have
\begin{equation}  \label{Eq7.3}
(q_{3}q_{1}^{\prime }q_{2}^{\prime n}=(q_{2}q_{1}^{\prime }q_{2}^{\prime
n}+(q_{1}q_{2}^{\prime }q_{3}^{\prime n}\ ,
\end{equation}
contradicting FLT, already proved true.

\subsection{Negative integers}

\label{Sc7.2}

We suppose all or some of the integers negative, and, \emph{at least two of
them, equal}, then through preliminary considerations, \emph{mutatis mutandis%
}, it is easy to conclude that, with no loss of generality, we still can
assume
\begin{equation}  \label{Eq7.4}
|z|>|x|>|y|\ .
\end{equation}

\subsubsection{All negative}

\label{Sc7.2.1}

If $n$ even, the situation is equivalent to FLT.

If $n$ odd, multiplying by $-1$, we have the same situation.

\subsubsection{Two negative}

\label{Sc7.2.2}

\begin{itemize}
\item[a)] $z$ and $x$
\begin{equation}  \label{Eq7.5}
(-z)^{n}=(-x)^{n}+y^{n}\ .
\end{equation}
If $n$ even, FLT.

If $n$ odd, the second number is a negative number with module smaller than $%
|z|$. Then \eqref{Eq7.5} is impossible. \vskip15pt

\item[b)] $z$ and $y$
\begin{equation}  \label{Eq7.6}
(-z)^{n}=x^{n}+(-y)^{n}\ .
\end{equation}
If $n$ even, FLT.

If $n$ odd, we have at right a positive number and at left a negative one,
an impossibility. \vskip15pt

\item[c)] $x$ and $y$
\begin{equation}  \label{Eq7.7}
z^{n}=(-x)^{n}+(-y)^{n}\ .
\end{equation}
If $n$ even, FLT.

If $n$ odd, positive at right, negative at left. Impossible.
\end{itemize}

\subsubsection{One negative}

\label{Sc7.2.3}

\begin{itemize}
\item[a)] $z$
\begin{equation}  \label{Eq7.8}
(-z)^{n}=x^{n}+y^{n}\ .
\end{equation}
If $n$ even, FLT.

If $n$ odd, negative at left, positive at right. Impossible. \vskip15pt

\item[b)] $x$
\begin{equation}  \label{Eq7.9}
z^{n}=(-x)^{n}+y^{n}\ .
\end{equation}
If $n$ even, FLT.

If $n$ odd, positive at left, negative at right. Impossible. \vskip15pt

\item[c)] $y$
\begin{equation}  \label{Eq7.10}
z^{n}=x^{n}+(-y)^{n}\ .
\end{equation}
If $n$ even, FLT.

If $n$ odd,
\begin{equation}  \label{Eq7.11}
x^{n}+(-y)^{n}<x^{n}<z^{n}\ ,
\end{equation}
thus, \eqref{Eq7.10} equality is impossible.
\end{itemize}

\subsection{Negative rationals}

\label{Sc7.3}

It is always possible to suppose that the signal of the rationals is that of
the $q_{i}^{\prime }$ in \eqref{Eq7.1}. Then it is easy to see that, after
multiplying \eqref{Eq7.1} by $q_{1}^{\prime }q_{2}^{\prime }q_{3}^{\prime }$%
, \eqref{Eq7.2} would contain, at least, one negative integer solution, an
impossibility, as already proved in the previous section.

\subsection{Algebraic irrationals and the inherent complex solutions}

\label{Sc7.4}

These exist and can be calculated through index $q>2$ roots of $y$, $x$ and $%
z$, in triplets admitting an equality with exponents $n\le 2$.

The first case is
\begin{equation}  \label{Eq7.12}
z=x+y\ .
\end{equation}

The candidate to be the solution is
\begin{equation}  \label{Eq7.13}
\{y^{1/q},x^{1/q},z^{1/q}\}\ .
\end{equation}

To begin with, we must know the relation between the $z$ root and the sum of
the other two. It cannot be
\begin{equation}  \label{Eq7.14}
z^{1/q}=x^{1/q}+y^{1/q}\ ,
\end{equation}
because this would imply the equalizing exponent being $1/q$, contradicting %
\eqref{Eq7.12}.

The exponent's increasing mechanism is $1/q,2/q,...$ up to $q/q=1$, giving
the equality \eqref{Eq7.12}.

Now, from Section~\ref{sec:3.3} we know that, once reversed with a certain
exponent, \emph{the inequality remains so with any greater exponents}.
Hence, if we admit
\begin{equation}  \label{Eq7.15}
z^{1/q}>x^{1/q}+y^{1/q}\ ,
\end{equation}
this inequality \emph{would remain reversed with the greater exponent $q/q=1$%
, contradicting \eqref{Eq7.12}}.

Therefore, the only remaining possibility is
\begin{equation}  \label{Eq7.16}
z^{1/q}<x^{1/q}+y^{1/q}\ ,
\end{equation}
so proving that \eqref{Eq7.13} is a solution with exponents $q>2$.

An example is $\{2^{1/3},3^{1/3},5^{1/3}\}$, which, powered to $3$, gives $%
2+3=5$, in this example with two inherent complex solutions.

The second case is the Pythagorean triplets
\begin{equation}
z^{2}=x^{2}+y^{2}\ .  \label{Eq7.17}
\end{equation}%
In total parallel terms,
\begin{equation}
\{y^{1/q},x^{1/q},z^{1/q}\}\ \ \ \ \hbox{and}\ \ \ \ q>1  \label{Eq7.18}
\end{equation}%
and the inherent complexes are solutions, giving \eqref{Eq7.17}, when
powered to $2q>2$.

\section{Numerical verification}

\label{sec:6}

\subsection{Proof}

\label{sec:6.1}

\begin{exmp}\label{Ex6.1}
\ $\{4,5,6\}$.
\begin{gather*}
 6<5+4;\ \ 6^{2}<5^{2}+4^{2};\ \ 6^{3}>5^{3}+4^{3} \ \ \to\ \ n=3;\ \ n-1=2\ponto\\
 T_{n-1}=T_{2}=\{4^{2},5^{2},6^{2}\}\ponto\\
\varphi=\frac{5^{2}+4^{2}}{6^{2}}=\frac{41}{36}=1.1388...\ponto\\
 k_{2}=\frac{5^{3}+4^{3}}{5^{2}+4^{2}}=\frac{189}{41}\ponto\\
 \frac{z}{k_{2}}=\frac{6\timesd41}{189}=1.3015...>\varphi=1.1388...\ponto\\
 \{\lambda\}=\left[\varphi,\frac{z}{k_{2}}\right]=\bigl[1.1388...,\,1.3015...\bigr]\ponto\\
 \frac{z}{k_{2}}\,z^{2}=\frac{6\timesd41\timesd36}{189}=46.8571...>5^{2}+4^{2}=41\ \ \to\\
 \to\ \ z^{3}>k_{2}\,(x^{2}+y^{2})=x^{3}+y^{3}\\
 6^{3}=216>5^{3}+4^{3}=189\ponto
\end{gather*}
\end{exmp}

\begin{exmp}\label{Ex6.2}
\ $\{8,9,10\}$.
\begin{gather*}
 10<9+8;\ \ 10^{2}<9^{2}+8^{2};\ \ ...;\ \ 10^{4}<9^{4}+8^{4};\ \ 10^{5}>9^{5}+8^{5}\ \ \to\\
 \to\ \ n=5;\ \ n-1=4\ponto\\
 T_{n-1}=T_{4}=\{8^{4},9^{4},10^{4}\}\ponto\\
 x^{n-1}+y^{n-1}=9^{4}+8^{4}=10\,657\ponto\\
 \varphi=\frac{x^{n-1}+y^{n-1}}{z^{n-1}}=\frac{9^{4}+8^{4}}{10^{4}}=\frac{10\,657}{10\,000}=1.0657\ponto
   \end{gather*}
 \begin{gather*}
 k_{4}=\frac{9^{5}+8^{5}}{9^{4}+8^{4}}=\frac{91\,817}{10\,657}\ponto\\
 \frac{z}{k_{4}}=10\timesd\frac{10\,657}{91\,817}=1.160678...>\varphi=1.0657\ponto\\
 \{\lambda\}=\left[\varphi,\frac{z}{k_{4}}\right]=\bigl[1.0657,\,1.160678...\bigr]\ponto\\
 \frac{z}{k_{4}}\,z^{4}=11\,606.78...>9^{4}+8^{4}=10\,657\ \ \to\\
 \to\ \ z^{5}>k_{4}\,(x^{4}+y^{4})\virg\\
 10^{5}>9^{5}+8^{5}\ponto
\end{gather*}
\end{exmp}

\subsection{Corollary 1}

\label{sec:6.2}

In each $\{y,x,z\}$ example it is indicated the set it belongs in Table~\ref%
{tab:1}.

\begin{exmp}\label{Ex6.3}
\ $\{2,5,9\}$.
\begin{gather*}
 9>5+2\ \ \to\ \ \mathrm{Set\ 1.1}\ \ \to\ \ n=1\pvirg \\
 b=\frac{\log(5+2)}{\log 9}=0.885...<1\ \ \to\ \ n=1\pvirg \\
 a=\frac{\log(5^{0}+2^{0})}{\log 9}=0.315...\pvirg \ \ \ \ b-a=0.885...-0.315...=0.570...\pvirg \\
 \frac{1}{2}<0.570...<1\ponto
\end{gather*}
\end{exmp}

\begin{exmp}\label{Ex6.4}
\ $\{2,7,9\}$.
\begin{gather*}
 9=7+2\ \ \to\ \ \mathrm{Set\ 1.2}\pvirg \\
 9^{2}>7^{2}+2^{2}\ \ \to\ \ n=2\pvirg \\
 b=\frac{\log(7^{2}+2^{2})}{\log 9}=1.806...\ \ \to\ \ n=2\pvirg \\
 a=\frac{\log(7+2)}{\log 9}=1\pvirg \ \ \ \ b-a=1.806...-1=0.806...\pvirg \\
 \frac{1}{2}<0.806...<1\ponto
\end{gather*}
\end{exmp}

\begin{exmp}\label{Ex6.5}
\ $\{4,5,7\}$.
\begin{gather*}
 7<5+4,\ 7^{2}>5^{2}+4^{2}\ \ \to\ \ \mathrm{Set\ 2.1}\ \ \to\ \ n=2\pvirg \\
 b=\frac{\log(5^{2}+4^{2})}{\log 7}=1.908...\ \ \to\ \ n=2\pvirg \\
 a=\frac{\log(5+4)}{\log 7}=1.129...\pvirg \ \ \ \ b-a=1.908...-1.129...=0.779...\pvirg \\
 \frac{1}{2}<0.779...<1\ponto
\end{gather*}
\end{exmp}

\begin{exmp}\label{Ex6.6}
\ $\{3,4,5\}$.
\begin{gather*}
 5<4+3,\ 5^{2}=4^{2}+3^{2},\ 5^{3}>4^{3}+3^{3}\ \ \to\ \ \mathrm{Set\ 2.2}\ \ \to\ \ n=3\pvirg \\
 b=\frac{\log(4^{3}+3^{3})}{\log 5}=2.802...\ \ \to\ \ n=3\pvirg \\
 a=\frac{\log(4^{2}+3^{2})}{\log 5}=2\pvirg \ \ \ \ b-a=2.802...-2=0.802...\pvirg \\
 \frac{1}{2}<0.802...<1\ponto
\end{gather*}
\end{exmp}

\begin{exmp}\label{Ex6.7}
\ $\{4,5,6\}$.
\begin{gather*}
6^{2}<5^{2}+4^{2},\ 6^{3}>5^{3}+4^{3}\ \ \to\ \ n=3\ \ \to\ \ \mathrm{Set\ 2.3.1}\pvirg \\
b=\frac{\log(5^{3}+4^{3})}{\log 6}=2.925...\ \ \to\ \ n=3\pvirg \\
a=\frac{\log(5^{2}+4^{2})}{\log 6}=2.072...\pvirg \ \ \ \ b-a=2.925...-2.072...=0.852...\pvirg \\
\frac{1}{2}<0.852...<1\ponto
\end{gather*}
\end{exmp}

\begin{exmp}\label{Ex6.8}
\ $\{6,7,8\}$.
\begin{align*}
 8<7+6,\ 8^{2}<7^{2}+6^{2},\ 8^{3}<7^{3}+6^{3},\ 8^{4}>7^{4}+6^{4}\ \
 &\to&\ \ \mathrm{Set\ 2.3.1}\\
 &\to&\ \ n=4\pvirg
\end{align*}
\vspace*{-4ex}
\begin{gather*}
b=\frac{\log(7^{4}+6^{4})}{\log 8}=3.950...\ \ \to\ \ n=4\pvirg \\
 a=\frac{\log(7^{3}+6^{3})}{\log 8}=3.042...\pvirg \ \ \ \ b-a=3.950...-3.042...=0.908...\pvirg \\
 \frac{1}{2}<0.908...<1\ponto
\end{gather*}
\end{exmp}

\begin{exmp}\label{Ex6.9}
\ $\{2,4,4\}$.
\begin{gather*}
 \mathrm{Set\ 2.3.1}\ \ \to\ \ \forall\,j\ \to\ 4^{j}<4^{j}+2^{j}\ \ \to\ \ \nexists\,n\pvirg \\
 b=\frac{\log(4^{n}+2^{n})}{\log 4}> n\,\frac{\log 4}{\log 4}>n\ \ \to\ \ \nexists\,n\!:\, n>n\ponto
\end{gather*}
\end{exmp}

\begin{exmp}\label{Ex6.10}
\ $\{3,3,3\}$.
\begin{gather*}
  \mathrm{Set\ 2.3.2}\ \ \to\ \ \forall\,j\ \to\ 3^{j}<3^{j}+3^{j}\ \ \to\ \ \nexists\,n\pvirg \\
    b=\frac{\log(3^{n}+3^{n})}{\log 3}=\log2+n\ \ \to\ \ \nexists\,n\!:\, n>\log2+n\ponto
\end{gather*}
\end{exmp}

\end{document}